\documentclass[11pt]{article}
\pdfoutput = 1

\usepackage[english]{babel}

\usepackage{geometry}
\usepackage{amsmath,mathtools}
\numberwithin{equation}{section} 

\usepackage[margin = 1em]{subfig}

\usepackage{graphicx} \graphicspath{{./}{Figures/}}

\usepackage{amsthm}

\newtheoremstyle{italic}
{5pt}
{5pt}
{\itshape}
{}
{}
{}
{.5em}
{\bfseries{\thmname{#1}~\thmnumber{#2}.}\thmnote{~\textnormal{(#3)}}}

\newtheoremstyle{upright}
{5pt}
{5pt}
{\upshape}
{}
{\bfseries}
{}
{.5em}
{\bfseries{\thmname{#1}~\thmnumber{#2}.}\thmnote{~\textnormal{(\textit{#3}\textrm{)}}}}

\theoremstyle{italic}
\newtheorem{theorem}{Theorem}[section]

\theoremstyle{upright}

\newtheorem{remark}[theorem]{Remark}

\usepackage{dsfont}
\usepackage{amsfonts}

\usepackage{xparse}

\newcommand\blfootnote[1]{%
  \begingroup
  \renewcommand\thefootnote{}\footnote{#1}%
  \addtocounter{footnote}{-1}%
  \endgroup
}

\newcommand{\Nat}{\mathbb{N}} 

\newcommand{\E}[1]{\mathbb{E}[#1]} 
\newcommand{\Prob}[1]{\mathbb{P}(#1)} 

\newcommand{\bld}[1]{\mathbf{#1}} 
\NewDocumentCommand \eb { m g }{%
    \IfNoValueTF{#2}
        {\bld{e}_{#1}}
        {\bld{e}^{(#1)}_{#2}} 
    }%
\NewDocumentCommand \zerob { g }{%
    \IfNoValueTF{#1}
        {\bld{0}}
        {\bld{0}^{(#1)}} 
    }%

\newcommand{\la}{\lambda} 

\newcommand{\be}{\beta}

\newcommand{\g}[2]{g_{#1;#2}} 
\newcommand{\f}[2]{f_{#1;#2}} 

\newcommand{\ST}[1]{B_{#1}} 
\newcommand{\BP}[1]{BP_{#1}} 

\newcommand{\ib}[1]{\bld{i}^{(#1)}} 
\newcommand{\jb}[1]{\bld{j}^{(#1)}} 
\newcommand{\qbs}[1]{\bld{q}^{(#1)}} 

\title{Joint queue length distribution of multi-class, single server queues with preemptive priorities}
\author{Andrei Sleptchenko\footnotemark[1], Jori Selen\footnotemark[2] \footnotemark[3], Ivo Adan\footnotemark[2] \footnotemark[3], and Geert-Jan van Houtum\footnotemark[3] \footnotemark[4]}

\begin{document}

\maketitle%

\renewcommand{\thefootnote}{\fnsymbol{footnote}}%
\footnotetext[1]{Department of Mechanical \& Industrial Engineering, Qatar University}%
\footnotetext[2]{Department of Mechanical Engineering, Eindhoven University of Technology}%
\footnotetext[3]{Department of Mathematics and Computer Science, Eindhoven University of Technology}%
\footnotetext[4]{School of Industrial Engineering, Eindhoven University of Technology}%
\blfootnote{E-mail address: {\tt j.selen@tue.nl}}%
\renewcommand{\thefootnote}{\arabic{footnote}} \setcounter{footnote}{0}%


\begin{abstract}%
In this paper we analyze an $M/M/1$ queueing system with an arbitrary number of customer classes, with class-dependent exponential service rates and preemptive priorities between classes. The queuing system can be described by a multi-dimensional Markov process, where the coordinates keep track of the number of customers of each class in the system. Based on matrix-analytic techniques and probabilistic arguments we develop a recursive method for the exact determination of the equilibrium joint queue length distribution. The method is applied to a spare parts logistics problem to illustrate the effect of setting repair priorities on the performance of the system. We conclude by briefly indicating how the method can be extended to an $M/M/1$ queueing system with non-preemptive priorities between customer classes.
\end{abstract}%


\section{Introduction}%
\label{sec:introduction}%

We consider a single server queueing system shared by $N$ customer classes, numbered $1,\ldots,N$. The class index $n$ indicates the priority rank; class 1 has the lowest priority and class $N$ has the highest priority. The arrival process of class-$n$ customers is a Poisson process with rate $\la_n$. The service time of class-$n$ customers is exponentially distributed with rate $\mu_n$. This system can be described by a multi-dimensional Markov process on the state space $\Nat_0^N$, where the coordinates keep track of the number of customers of each class in the system.

In this paper we present an \textit{exact} method, based on matrix-analytic techniques \cite{neuts1989_matrix-analytic_approach,ramaswami1988_matrix-analytic_approach} to determine the equilibrium joint queue length distribution. In particular, it appears to be possible to avoid the use of infinite series and truncation of the state space. The crucial observation is that the Markov process, embedded on states in which there are no customers of priority classes higher than $n$, is of the $M/G/1$ type, where the number of class-$n$ customers represents the class-$n$ level. This is due to the fact that during excursions of the Markov process in which higher priority customers are present, any number of lower priority customers may arrive. Thus, a natural way to find the equilibrium joint queue length distribution is by recursive application of the theory of $M/G/1$-type Markov processes.

The joint queue length distribution is required in applications in the area of spare parts logistics and production. Specifically, our interest in the $M/M/1$ priority system with $N$ classes arose from a spare parts logistics problem, where the joint queue length distribution is necessary for an exact performance analysis. This problem is discussed in Section~\ref{sec:application}.

Priority queueing systems have a long history (cf. \cite{cobham1954_priority_assignment,davis1966_old_marginal,jaiswal1968_priority_queues}) and single and multi-server priority queues received much attention. Most of the earlier studies concentrate on the transforms of \textit{marginal} system characteristics such as the queue length and waiting time of a specific priority class. The focus on marginal system characteristics is also seen in recent work in \cite{horvath2005_approximate_G-G-1,vuuren2007_matrix-analytic-G-G-1_Nclass}, where the domain of priority queueing systems with general arrival and service time distributions is treated.

Joint queue length distributions have first been studied in \cite{miller1981_matrix-geometric_2class} using the matrix-geometric method \cite{neuts1981_matrix_geometric_approach} for an $M/M/1$ priority queueing system with two classes. This study spurred the observation made in \cite{alfa1998_matrix-geometric_method_2classes,wagner1996_matrix-geometric_MAP-PH-c,wagner1998_matrix-geometric_MAP-PH-c} that the matrix-geometric method is a natural choice for studying priority queueing systems with a quasi-birth--death (QBD) structure. In these papers, the matrix-geometric method is generalized to systems with two priority classes, a Markovian arrival process and a phase-type service time distribution. In \cite{alfa2003_matrix-geometric_method_Nclasses} the same matrix-geometric method is applied to a discrete-time $N$-class system, leading to an approximation of the joint equilibrium distribution, as the rate matrix $R$ needs to be truncated for actual computation. An $M/PH/1$ non-preemptive priority system with $N$ classes with different service rates per class is studied in \cite{isotupa2002_matrix-geometric_M-PH-1}, where an algorithm is derived using matrix-geometric techniques for the computation of the joint queue length distribution for three aggregated classes. The observation that is \textit{not} made in \cite{alfa2003_matrix-geometric_method_Nclasses,isotupa2002_matrix-geometric_M-PH-1} is that lower priority customers see the queueing system as an $M/G/1$-type system, i.e., an $M/M/1$ system with an unreliable server (or vacations), where down times correspond to high priority service interruptions. This observation is made and implemented in \cite{harchol2005_reduce_dimensions,wierman2006_dimensionality_reduction_M-PH-c}, where the distribution of the down times are approximated by phase-type distributions, the first three moments of which are matched to the moments of high priority service interruptions. However, only marginal queue length distributions are obtained.

There is also a number of papers studying the joint queue length distribution using alternative approaches. Generating functions are used in \cite{gail1988_GF_non-preemptive,gail1992_GF_preemptive} for the analysis of $M/M/c$ priority queueing systems with two classes. Generating functions are also used in \cite{mitrani1981_GF_M-M-c} for an $M/M/c$ preemptive priority system with more than two classes. Here, customers of higher priority are aggregated, leading to an approximation of the equilibrium distribution. Later, \cite{sleptchenko2003_nonpre_approximation_M-M-c_2class,sleptchenko2002_pre_approximation_M-M-c_2class,sleptchenko2005_GF_matrix-geometric_M-M-c_2class} use a mixture of the matrix-geometric method and generating function technique to analyze preemptive and non-preemptive priority $M/M/c$ queueing systems with two classes, where each class can have different types of customers. The mixture of the two methods leads to an approximation of the joint equilibrium distribution as the number of matrix operations has to be finite for actual computation.

The area of priority queueing systems still is an active field of research. More recently, priority queueing systems with impatient high priority customers have been analyzed using generating functions \cite{brandt2004_Laplace_GF_M-M-1+GI}; by identifying simple Markov processes \cite{choi2001_M-M-1+D}; using a level-crossing method \cite{iravani2008_level-crossing_M-GI-1+M} or using Laplace-Stieltjes transforms \cite{jouini2013_Laplace_Stieltjes_M-M-c+M}. These systems have applications in, for example, telecommunication systems where voice messages need to be delivered timely and have priority over data packets. An alternative to impatient customers are queueing systems where customers can reduce their sojourn time by transferring to a higher priority class. This allows impatient customers to be served earlier. In \cite{xie2009_stationary_distribution_transfers}, bounds on the equilibrium distribution are given. The study of a queueing system with transferring customers is motivated by the potential application in the design of emergency departments. Here, patients are categorized in classes of different priority, where patients can transfer from a lower priority class to a higher priority class. Approximations for the first and second moment of the waiting time in an $M/G/c$ non-preemptive priority queueing system with an arbitrary number of priority classes are given in \cite{alhanbali2013_approximations_waiting_time_distribution_M-G-c}.

Our main contribution is that we describe a method for the exact determination of the joint queue length distribution for a preemptive priority queueing system with an arbitrary number of classes and class-dependent service rates. We use the property that the embedded Markov process is of the $M/G/1$ type. Key to the approach is identifying first passage probabilities which are computed by one-step analysis. We then recursively apply matrix-analytic methods related to $M/G/1$-type Markov processes and avoid the use of infinite series.

The remainder of the paper is organized as follows. In Section~\ref{sec:matrix-analytic_method} we describe how the matrix-analytic method is applied to an $N$-class preemptive priority single server system. To ease the understanding of the method in general and highlight the recursive nature, we first treat the two and three-class systems in Sections~\ref{subsec:two-class_system} and \ref{subsec:three-class_system}, respectively. Next, in Section~\ref{sec:application} we present the application in spare parts logistics where the joint queue length distribution is needed for an exact analysis. In the final section we conclude by indicating how to extend the method to non-preemptive priority rules.


\section{Matrix-analytic method}%
\label{sec:matrix-analytic_method}%

The $M/M/1$ preemptive priority system can be described by a Markov process with states $(q_N,\ldots,q_1)$, where $q_n$ denotes the number of class-$n$ customers in the system. State transitions are triggered by arrival and service completions. Class-$n$ customers enter at an exponential rate $\la_n$, triggering a transition from $(q_N,\ldots,q_1)$ to state $(q_N,\ldots,q_n + 1,\ldots,q_1)$, and if $q_N = \cdots = q_{n+1} = 0$ and $q_n > 0$, class-$n$ customers are served at an exponential rate $\mu_n$, which leads to a transition from $(0,\ldots,0,q_n,\ldots,q_1)$ to $(0,\ldots,0,q_n - 1,\ldots,q_1)$. Throughout the paper we assume that the system is stable, i.e., the traffic intensity $\rho$ is less than 1 (see, e.g., \cite{gross1974_fundamentals}):
\begin{equation}%
\rho \coloneqq \sum_{i = n}^N \la_n / \mu_n < 1, \label{eqn:load_rho}
\end{equation}%
and we denote by $p(q_N,\ldots,q_1)$ the equilibrium probability of being in state $(q_N,\ldots,q_1)$. To ease notation, let us introduce $\la \coloneqq \sum_{n = 1}^N \la_n$. We propose to use the matrix-analytic method for $M/G/1$ structured systems to exactly and recursively calculate the joint queue length probabilities $p(q_N,\ldots,q_1)$, starting from $p(0,\ldots,0) = 1 - \rho$. Key to this approach are first passage probabilities, that can be determined through one-step analysis. In fact, the first passage probabilities are the elements of the auxiliary matrix $G$ of the matrix-analytic method. However, rather than determining the infinite matrix $G$ using matrix equations, we recursively determine its elements using scalar equations, derived by exploiting the skip-free property of this Markov process. To highlight the recursive nature of the method we first treat the two and three-class systems.


\subsection{Two-class system}%
\label{subsec:two-class_system}%

\begin{figure}%
\centering%
\subfloat[Transition rate diagram of the two-class system.]{%
\includegraphics{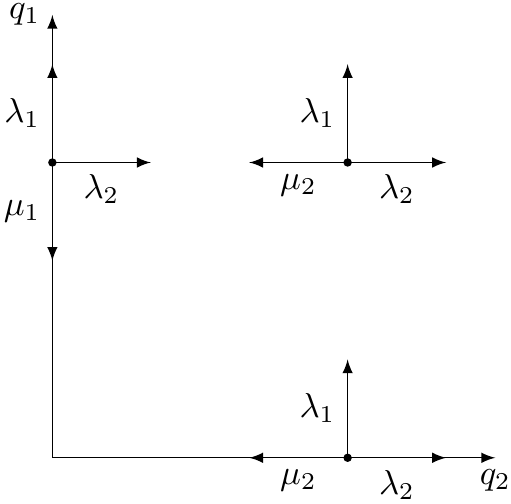}
}%
\subfloat[Embedded on class-2 level $q_2 = 0$.]{%
\includegraphics{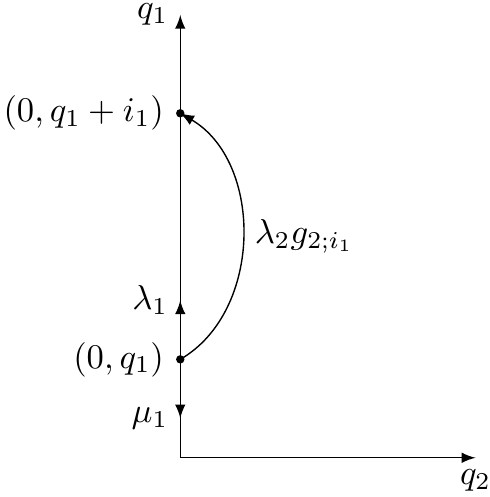}
}%
\caption{Transition rate diagrams of the two-class system.}
\label{fig:two-class_system}%
\end{figure}%

The transition rate diagram of the two-class system depicted in Figure~\ref{fig:two-class_system}(a) shows that the two-class system is a QBD process with class-2 levels $q_2$ defined as the set of states with $q_2$ high priority customers. To calculate the probabilities $p(q_2,q_1)$, we propose to exploit the $M/G/1$ structure of this Markov process, instead of its $G/M/1$ structure as done by \cite{miller1981_matrix-geometric_2class}. Instrumental in the calculation of $p(q_2,q_1)$ are the first passage probabilities $\g{2}{i_1}$, instead of the elements of the rate matrix as in \cite{miller1981_matrix-geometric_2class}. The first passage probability $\g{2}{i_1}$ is defined as the probability that, starting at class-2 level $q_2 > 0$ in state $(q_2,q_1)$, the first passage to class-2 level $q_2 - 1$ happens in state $(q_2 - 1,q_1 + i_1)$. Note that $\g{2}{i_1}$ does not depend on the starting state $(q_2,q_1)$, and can be interpreted as the probability that $i_1$ class-1 customers arrive during a busy period of class-2 customers. By one-step analysis we get
\begin{align}%
\mu_2 - (\la + \mu_2)\g{2}{0} + \la_2 \g{2}{0}^2 &= 0, \quad i_1 = 0, \label{eqn:two-class_g_2;0} \\
- (\la + \mu_2)\g{2}{i_1} + \la_1 \g{2}{i_1 - 1} + \la_2 \sum_{j_1 = 0}^{i_1} \g{2}{j_1} \g{2}{i_1-j_1} &= 0, \quad i_1 > 0. \label{eqn:two-class_g_2;i_1}
\end{align}%
So $\g{2}{i_1}$ can be recursively calculated, starting from $\g{2}{0}$, which follows from \eqref{eqn:two-class_g_2;0},
\begin{equation}%
\g{2}{0} = \frac{1}{2\la_2} \Bigl( \la + \mu_2 - \bigl((\la + \mu_2)^2 - 4\la_2\mu_2 \bigr)^{\frac{1}{2}} \Bigr). \label{eqn:two-class_g_2;0_explicit}
\end{equation}%
To calculate $p(q_2,q_1)$, we use the following equation for excursions starting at class-2 level $q_2$ to levels higher than $q_2$ ending at first return to class-2 level $q_2$. The number of excursions per time unit that end in state $(q_2,q_1)$ is equal to $p(q_2 + 1,q_1)\mu_2$, but this number is also equal to the excursions starting from class-2 level $q_2$ per time unit that end in state $(q_2,q_1)$. The number of excursions per time unit that start in state $(q_2,q_1 - i_1)$ is equal to $p(q_2,q_1 - i_1)\la_2$, a fraction $\g{2}{i_1}$ of which ends in $(q_2,q_1)$. Hence,
\begin{equation}%
p(q_2 + 1,q_1)\mu_2 = \sum_{i_1 = 0}^{q_1} p(q_2,q_1 - i_1) \la_2 \g{2}{i_1}, \quad q_2,q_1 \ge 0, \label{eqn:two-class_p(q_2,q_1)}
\end{equation}%
from which all probabilities can be recursively calculated, once the boundary probabilities $p(0,q_1)$ are known. The probabilities $p(0,q_1)$ can be determined by considering the Markov process embedded on class-2 level 0. The transition rate diagram of the embedded Markov process is shown in Figure~\ref{fig:two-class_system}(b). Note that the embedded Markov process has an $M/G/1$ structure with class-1 levels $q_1$ defined as the set of states with $q_1$ class-1 customers (and no class-2 customers). To formulate the analogue of \eqref{eqn:two-class_p(q_2,q_1)}, we introduce $\f{2}{i_1}$, which is the probability that, starting in state $(1,q_1)$, the first passage to class-1 levels less than or equal to $q_1 + i_1$ happens in state $(0,q_1 + i_1)$. In this case, this first passage probability is equal to the probability that during a busy period of class-2 customers, at least $i_1$ class-1 customers arrive. So
\begin{equation}%
\f{2}{i_1} = 1 - \sum_{j_1 = 0}^{i_1 - 1} \g{2}{j_1}. \label{eqn:two-class_f_2;i_1}
\end{equation}%
Then, similar to \eqref{eqn:two-class_p(q_2,q_1)}, we have
\begin{equation}%
p(0,q_1 + 1)\mu_1 = p(0,q_1)\la_1 + \sum_{i_1 = 0}^{q_1} p(0,q_1 - i_1) \la_2 \f{2}{i_1 + 1}, \quad q_1 \ge 0, \label{eqn:two-class_p(0,q_1)}
\end{equation}%
which can be used to calculate all boundary probabilities, starting from the probability of an empty system $p(0,0) = 1 - \rho$.


\subsection{Three-class system}%
\label{subsec:three-class_system}%

\begin{figure}%
\centering%
\subfloat[$q_3 > 0$.]{%
\includegraphics{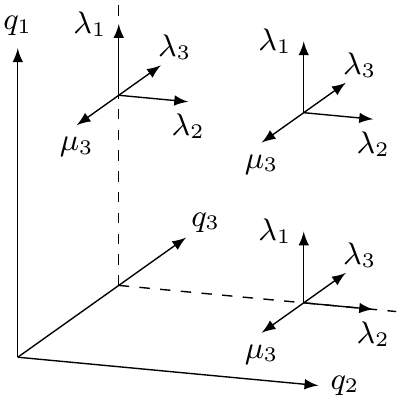}
}%
\subfloat[$q_3 = 0$.]{%
\includegraphics{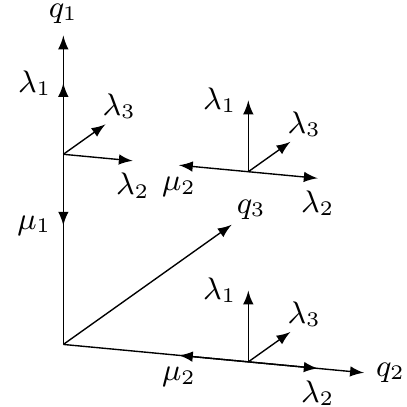}
}%
\caption{Transition rate diagram of the three-class system.}%
\label{fig:three-class_system}
\end{figure}%

The transition rate diagram of the three-class system is shown in Figure~\ref{fig:three-class_system}(a) and (b). This system can be described by a QBD process with class-3 levels $q_3$ defined as the set of states with $q_3$ high priority customers. Let $\g{3}{i_2,i_1}$ be the probability that, starting at class-3 level $q_3 > 0$ in state $(q_3,q_2,q_1)$, the first passage to class-3 level $q_3 - 1$ happens in state $(q_3 - 1,q_2 + i_2, q_1 + i_1)$. Note that $\g{3}{i_2,i_1}$ can be interpreted as the probability that $i_2$ class-2 and $i_1$ class-1 customers arrive during a busy period of high priority class-3 customers. By one-step analysis,
\begin{align}%
\mu_3 - ( \la + \mu_3 ) \g{3}{0,0} + \la_3 \g{3}{0,0}^2 &= 0, \quad i_2,i_1 = 0, \label{eqn:three-class_g_3;0,0} \\
- ( \la + \mu_3 )\g{3}{i_2,i_1} + \la_1 \g{3}{i_2,i_1 - 1} + \la_2 \g{3}{i_2 - 1,i_1} & \notag \\
+ \la_3 \sum_{j_2 = 0}^{i_2} \sum_{j_1 = 0}^{i_1} \g{3}{j_2,j_1} \g{3}{i_2 - j_2,i_1 - j_1} &= 0, \quad i_2 + i_1 > 0, \label{eqn:three-class_g_3;i_2,i_1}
\end{align}%
where by convention, $\g{3}{i_2,i_1} = 0$ if $i_2 < 0$ or $i_1 < 0$. From \eqref{eqn:three-class_g_3;i_2,i_1} the probabilities $\g{3}{i_2,i_1}$ can be recursively calculated, starting from $\g{3}{0,0}$, which follows from \eqref{eqn:three-class_g_3;0,0},
\begin{equation}%
\g{3}{0,0} = \frac{1}{2\la_3} \Bigl( \la + \mu_3 + \bigl( ( \la + \mu_3 )^2 - 4 \la_3 \mu_3 \bigr)^{\frac{1}{2}} \Bigr). \label{eqn:three-class_g_3;0,0_explicit}
\end{equation}%
Similar to \eqref{eqn:two-class_p(q_2,q_1)}, we have
\begin{equation}%
p(q_3 + 1,q_2,q_1)\mu_3 = \sum_{i_2 = 0}^{q_2} \sum_{i_1 = 0}^{q_1} p(q_3,q_2 - j_2,q_1 - j_1) \la_3 \g{3}{i_2,i_1}, \quad q_3,q_2,q_1 \ge 0, \label{eqn:three-class_p(q_3,q_2,q_1)}
\end{equation}%
which can be utilized to calculate all probabilities, once the boundary probabilities $p(0,q_2,q_1)$ are known.

\begin{figure}%
\centering%
\subfloat[Embedded on class-3 level 0.]{%
\includegraphics{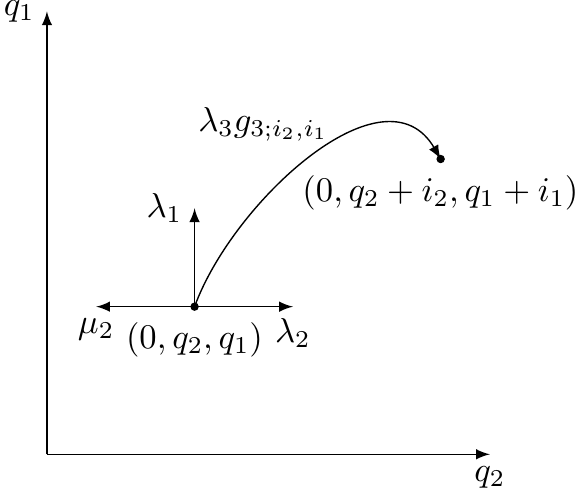}
}%
\subfloat[Embedded on the axis $q_3 = q_2 = 0$.]{%
\includegraphics{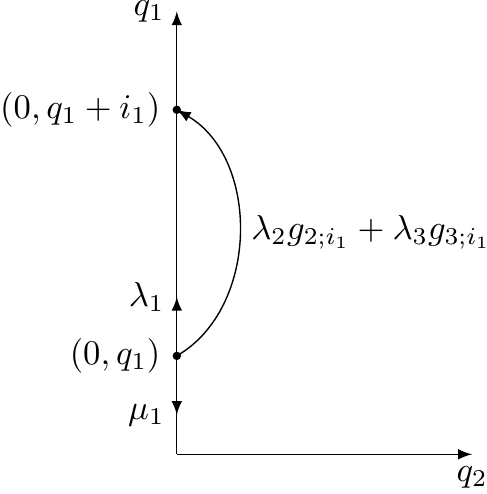}
}%
\caption{Transition rate diagram of two embedded Markov processes of the three-class system.}%
\label{fig:three-class_system_embedded}%
\end{figure}%

To determine $p(0,q_2,q_1)$ we proceed by considering the Markov process embedded on class-3 level 0, which is of the $M/G/1$ type, with class-2 levels $q_2$ defined as the set of states with $q_2$ class-2 customers (and no class-3 customers). Its transition rate diagram is depicted in Figure~\ref{fig:three-class_system_embedded}(a). The first passage probabilities $\g{2}{i_1}$ for the embedded Markov process are defined as the probability that, when starting at class-2 level $q_2 > 0$ in state $(0,q_2,q_1)$, the first passage to class-2 level $q_2 - 1$ happens in state $(0,q_2 - 1,q_1 + i_1)$. Further, the first passage probabilities $\g{3}{i_1}$ are defined as the probability that, when starting in state $(1,q_2 - 1,q_1)$, the first passage to class-2 level $q_2 - 1$ happens in state $(0,q_2 - 1,q_1 + i_1)$. Observe that $\g{k}{i_1}$ is the probability that $i_1$ class-1 customers arrive during a busy period of higher priority (class-2 and class-3) customers, that starts with the arrival of a class-$k$ customer, for $k = 2,3$. Notice the difference between the first passage probabilities $\g{3}{i_1}$ and $\g{3}{i_2,i_1}$. The number of indices after the semicolon in the subscript is related to what level the Markov process is embedded on, as illustrated in Figure~\ref{fig:three-class_system_embedded}. By one-step analysis we get for $k = 2,3$,
\begin{align}%
\mu_k - ( \la + \mu_k ) \g{k}{0} + \sum_{m = 2}^3 \la_m \g{m}{0} \g{k}{0} &= 0, \quad i_1 = 0, \label{eqn:three-class_g_k;0} \\
- ( \la + \mu_k ) \g{k}{i_1} + \la_1 \g{k}{i_1-1} + \sum_{m = 2}^3 \la_m \sum_{j_1 = 0}^{i_1} \g{m}{j_1} \g{k}{i_1 - j_1} &= 0, \quad i_1 > 0. \label{eqn:three-class_g_k;i_1}
\end{align}%
From equations \eqref{eqn:three-class_g_k;i_1}, both $\g{2}{i_1}$ and $\g{3}{i_1}$ can be recursively calculated, with $\g{2}{0}$ and $\g{3}{0}$ being the minimal nonnegative solution of \eqref{eqn:three-class_g_k;0}. To solve \eqref{eqn:three-class_g_k;0} we introduce $\ST{k}$ which is the Laplace-Stieltjes transform (LST) of the service time of a class-$k$ customer and $\BP{k}$ the LST of a high priority (class-2 and class-3) busy period initiated by a class-$k$ customer. Then $\g{2}{0}$ and $\g{3}{0}$ can be calculated from, see \cite[Section~5.8]{kleinrock1975_queueing_systems},
\begin{align}%
\g{k}{0} &= \BP{k}(\la_1) = \ST{k}(\la_1 + (\la_2 + \la_3)(1 - \BP{2,3}(\la_1))) \notag \\
&= \frac{\mu_k}{\mu_k + \la_1 + (\la_2 + \la_3)(1 - \BP{2,3}(\la_1))}, \quad k = 2,3, \label{eqn:three-class_g_k;0_LST}
\end{align}%
where $\BP{2,3}(s)$ is the LST of a high priority busy period initiated by a class-2 or a class-3 customer, which is equal to the LST of the busy period in an $M/H_2/1$ queue with class-2,3 customers,
\begin{equation}%
\BP{2,3}(s) = \sum_{m = 2}^3 \frac{\la_m}{\la_2 + \la_3} \frac{\mu_m}{\mu_m + s + (\la_2 + \la_3)(1 - \BP{2,3}(s))}, \quad s \ge 0. \label{eqn:three-class_BP_2,3}
\end{equation}%
To formulate the analogue of \eqref{eqn:three-class_p(q_3,q_2,q_1)} for $p(0,q_2,q_1)$, we introduce the first passage probabilities $\f{3}{i_2,i_1}$ defined as the probability that, when starting in state $(1,q_2,q_1)$, the first passage to class-2 levels less than or equal to $q_2 + i_2$ happens in state $(0,q_2 + i_2,q_1 + i_1)$. The probability $\f{3}{i_2,i_1}$ can be interpreted as the probability that at the end of a busy period of class-3 customers, there have been at least $i_2$ class-2 arrivals, and then, when the server brings down this number to $i_2$, the total number of class-1 arrivals (from the start of the busy period of class-3 customers) has been $i_1$. Hence, we can express $\f{3}{i_2,i_1}$ as an infinite sum,
\begin{equation}%
\f{3}{i_2,i_1} = \sum_{m = 0}^\infty \sum_{\substack{j_0,\ldots,j_m \ge 0 \\ j_0 + \cdots + j_m = i_1}} \g{3}{i_2 + m,j_0} \g{2}{j_1} \cdots \g{2}{j_m}. \label{eqn:three-class_f_3;i_2,i_1_definition}
\end{equation}%
Before elaborating on the computation of $\f{3}{i_2,i_1}$, we proceed to derive an equation for the probabilities $p(0,q_2,q_1)$, by considering excursions to class-2 levels higher than $q_2$, that start at class-2 level $q_2$ or lower, and end at first return to class-2 level $q_2$ in state $(0,q_2,q_1)$. The number of excursions per time unit that end in state $(0,q_2,q_1)$ is equal to $p(0,q_2 + 1,q_1)\mu_2$. This number is also equal to the excursions starting from class-2 level $q_2$ or lower per time unit that end in state $(0,q_2,q_1)$. A fraction $\g{2}{i_1}$ of the excursions starting in $(0,q_2,q_1 - i_1)$ by a class-2 arrival end in $(0,q_2,q_1)$. Excursions to class-2 levels higher than $q_2$ starting in state $(0,q_2 - i_2,q_1 - i_1)$ by a class-3 arrival reach, with probability $\f{3}{i_2,i_1} - \g{3}{i_2,i_1}$, class-2 level $q_2$ in state $(0,q_2,q_1)$ at first return to class-2 level $q_2$. Note that $\g{3}{i_2,i_1}$ needs to be subtracted, since with probability $\g{3}{i_2,i_1}$ class-2 level $q_2$ is reached, but not yet exceeded. Hence,
\begin{align}%
p(0,q_2 + 1,q_1)\mu_2 &= \sum_{i_1 = 0}^{q_1} p(0,q_2,q_1 - i_1) \la_2 \g{2}{i_1} \notag \\
&\quad + \sum_{i_2 = 0}^{q_2} \sum_{i_1 = 0}^{q_1} p(0,q_2 - i_2,q_1 - i_1) \la_3 (\f{3}{i_2,i_1} - \g{3}{i_2,i_1}), \quad q_2,q_1 \ge 0, \label{eqn:three-class_p(0,q_2,q_1)}
\end{align}%
from which $p(0,q_2,q_1)$ can be recursively calculated, once the boundary probabilities $p(0,0,q_1)$ are known. To determine $p(0,0,q_1)$ we consider the Markov process embedded on the axis $q_3 = q_2 = 0$, the transition rate diagram of which is depicted in Figure~\ref{fig:three-class_system_embedded}(b), with class-1 levels $q_1$ defined as the set of states with $q_1$ class-1 customers (and no class-2 or class-3 customers). To finally formulate the equations for $p(0,0,q_1)$ we define $\f{k}{i_1}$ as the probability that, when starting in state $(0,1,q_1)$ if $k = 2$ and starting in state $(1,0,q_1)$ if $k = 3$, the first passage to class-1 levels less than or equal to $q_1 + i_1$ happens in state $(0,0,q_1 + i_1)$. Similar as for the two-class system, this first passage probability is equal to the probability that at least $i_1$ class-1 customers arrive during a busy period of class-2,3 customers, initiated by a class-$k$ customer. So, for $k = 2,3$,
\begin{equation}%
\f{k}{i_1} = 1 - \sum_{j_1 = 0}^{i_1 - 1} \g{k}{j_1}. \label{eqn:three-class_f_k;i_1}
\end{equation}%
Then, similar to \eqref{eqn:two-class_p(0,q_1)}, we have
\begin{equation}%
p(0,0,q_1 + 1)\mu_1 = p(0,0,q_1)\la_1 + \sum_{i_1 = 0}^{q_1} p(0,0,q_1 - i_1)(\la_2 \f{2}{i_1 + 1} + \la_3 \f{3}{i_1 + 1}), \quad q_1 \ge 0. \label{eqn:three-class_p(0,0,q_1)}
\end{equation}%
This equation can be used to recursively calculate $p(0,0,q_1)$, with initially $p(0,0,0) = 1 - \rho$.

We now turn to the calculation of the first passage probabilities $\f{3}{i_2,i_1}$. To avoid evaluation of the infinite sums in \eqref{eqn:three-class_f_3;i_2,i_1_definition}, we again employ one-step analysis, yielding for $i_2 > 0$ and $i_1 \ge 0$,
\begin{align}%
-( \la + \mu_3 ) \f{3}{i_2,i_1} + \la_1 \f{3}{i_2,i_1 - 1} + \la_2 \f{3}{i_2 - 1,i_1} & \notag \\
+ \la_3 \Bigl( \sum_{j_2 = 0}^{i_2 - 1} \sum_{j_1 = 0}^{i_1} \g{3}{j_2,j_1} \f{3}{i_2 - j_2,i_1 - j_1}  + \sum_{j_1 = 0}^{i_1} \f{3}{i_2,j_1} \g{3}{i_1 - j_1} \Bigr) &= 0, \label{eqn:three-class_f_3;i_2,i_1}
\end{align}%
where by convention, $\f{3}{i_2,i_1} = 0$ if $i_1 < 0$. The first passage probabilities $f_{3;i_2,i_1}$ can be recursively calculated using the equations \eqref{eqn:three-class_f_3;i_2,i_1}, starting with $\f{3}{0,i_1} = \g{3}{i_1}$. The last two terms in \eqref{eqn:three-class_f_3;i_2,i_1} need some explanation: this is the probability of first passage to class-2 levels less than or equal to $q_2 + i_2$ in state $(0,q_2 + i_2,q_1 + i_1)$ when starting an excursion in state $(2,q_2,q_1)$, so with two instead of one class-3 customer. Now imagine that the second class-3 customer enters service when the busy period generated by the first class-3 customer finishes. The first term corresponds to the event that the number of class-2 arrivals during the busy period generated by the first class-3 customer is $j_2 < i_2$, so that the number of class-2 arrivals during the second busy period should be at least $i_2 - j_2$. The second term corresponds to the event that the number of class-2 arrivals during the first busy period is $j_2 \ge i_2$. The surplus number $j_2 - i_2$ of class-2 customers should be served after the busy period generated by the second class-3 customer. The duration of the excursion will not be altered if these class-2 customers enter service (as well as any higher priority customer arriving during their service) before the second class-3 customer. Then $\f{3}{i_2,j_1}$ is the probability that the number of class-1 arrivals is $j_1$ when the last surplus class-2 customer completes service, and thus the number of class-1 arrivals during the busy period generated by the second class-3 customer should be exactly equal to $i_1 - j_1$. Note that the busy period generated by this second class-3 customer includes class-3 and class-2 customers, since each arriving class-2 customer is surplus. So the probability of exactly $i_1 - j_1$ class-1 arrivals is $\g{3}{i_1 - j_1}$.


\subsection{$N$-class system}%
\label{subsec:N-class_system}%

We now extend the approach for obtaining the stationary distribution of the three-class system to an $N$-class system. Since we are dealing with $N$ classes, we need some accommodating notation. We introduce $\ib{n} = (i_n,i_{n - 1},\ldots,i_1)$, $\qbs{n} = (q_n,q_{n - 1},\ldots,q_1)$, $\jb{n}$ is vector-index of length $n$, $\zerob{n}$ is the zero vector of length $n$ and $\eb{n}{k}$ denotes a vector of zeros of length $n$ with a 1 at position $n + 1 - k$. Class-$n$ level $q_n$ denotes the set of states with $q_n$ class-$n$ customers and no customers of higher classes.

Once again, we have two types of first passage probabilities. The first type is the first passage probability $\g{k}{\ib{n}}, ~ k \ge n + 1$ defined as the probability that, when starting in state $(0,\ldots,0,q_{n + 1} - 1,\qbs{n}) + \eb{N}{k}$, the first passage to class-$(n + 1)$ level $q_{n + 1} - 1$ happens in state $(0,\ldots,0,q_{n + 1} - 1,\qbs{n} + \ib{n})$. Second, $\f{k}{\ib{n}}, ~ k \ge n + 1$ is the probability that, when starting in state $(0,\ldots,0,\qbs{n}) + \eb{N}{k}$, the first passage to class-$n$ levels less than or equal to $q_{n} + i_{n}$ happens in state $(0,\ldots,0,\qbs{n} + \ib{n})$.

We first describe how to obtain the first passage probabilities, followed by the computation of the equilibrium probabilities. Note that $\f{k}{0,\ib{n - 1}} = \g{k}{\ib{n - 1}}$. By one-step analysis we get for $n = N - 1,N - 2,\ldots,1$ and $k \ge n + 1$,
\begin{align}%
\mu_k - (\la + \mu_k)\g{k}{\zerob{n}} + \sum_{m = n + 1}^N \la_m \g{m}{\zerob{n}} \g{k}{\zerob{n}} &= 0, \quad \ib{n} = \zerob{n}, \label{eqn:N-class_g_k;0} \\
-(\la + \mu_k)\g{k}{\ib{n}} + \sum_{m = 1}^{n} \la_m \g{k}{\ib{n} - \eb{n}{m}}& \notag \\
+ \sum_{m = n + 1}^N \la_m \sum_{\jb{n} = \zerob{n}}^{\ib{n}} \g{m}{\jb{n}} \g{k}{\ib{n} - \jb{n}} &= 0, \quad \sum_{m = 1}^n i_m > 0. \label{eqn:N-class_g_k;i(n-1)}
\end{align}%
From \eqref{eqn:N-class_g_k;i(n-1)}, all $\g{k}{\ib{n}}$ with $k \ge n + 1$ and $n$ fixed can be calculated, with $\g{k}{\zerob{n}}$ computed as
\begin{align}%
\g{k}{\zerob{n}} &= \BP{k}\Bigl(\sum_{m = 1}^{n} \la_m \Bigr) = \ST{k}\Bigl(\sum_{m = 1}^{n} \la_m + \sum_{m = n + 1}^N \la_m (1 - \BP{n + 1,\ldots,N}(\sum_{m = 1}^{n} \la_m )) \Bigr) \notag \\
&= \frac{\mu_k}{\mu_k + \sum_{m = 1}^{n} \la_m + \sum_{m = n + 1}^N \la_m (1 - \BP{n + 1,\ldots,N}(\sum_{m = 1}^{n} \la_m))}, \label{eqn:N-class_g_k;0_explicit}
\end{align}%
where $\BP{n + 1,\ldots,N}(s)$ is the LST of a high priority (class-$(n+1)$ and higher) busy period, which is equal to the LST of the busy period in an $M/H_{N-n}/1$ queue with class-$(n+1),\ldots,N$ customers,
\begin{align}%
\BP{n + 1,\ldots,N}(s) &= \sum_{m = n + 1}^{N} \frac{\la_m}{\sum_{l = n + 1}^N \la_l} \frac{\mu_m}{\mu_m + s + \sum_{l = n + 1}^N \la_l(1 - \BP{n + 1,\ldots,N}(s))}, \quad s \ge 0.
\end{align}%

The first passage probabilities $\f{k}{\ib{n}}$ with $n = N - 1,\ldots,2$ and $k \ge n + 1$ follow from one-step analysis similar to \eqref{eqn:three-class_f_3;i_2,i_1}, with $i_n > 0$,
\begin{align}%
-(\la + \mu_k)\f{k}{\ib{n}} &+ \sum_{m = 1}^{n} \la_m \f{k}{\ib{n} - \eb{n}{m}} + \sum_{m = n + 1}^N \la_m \Bigl( \sum_{j_{n} = 0}^{i_{n} - 1} \sum_{\jb{n - 1} = \zerob{n - 1}}^{\ib{n - 1}} \g{m}{j_n,\jb{n - 1}} \f{k}{i_n - j_n} \notag \\
&+ \sum_{\jb{n - 1} = \zerob{n - 1}}^{\ib{n - 1}} \f{m}{i_{n},\jb{n - 1}} \g{k}{\ib{n - 1} - \jb{n - 1}} \Bigr) = 0, \quad \ib{n - 1} \ge \zerob{n - 1}. \label{eqn:N-class_f_k;i(n-1)}
\end{align}%
The last two terms in \eqref{eqn:N-class_f_k;i(n-1)} describe the probability of first passage to class-$n$ levels less than or equal to $q_{n} + i_{n}$ in state $(\zerob{N - n},\qbs{n} + \ib{n})$ when starting an excursion in state $(\zerob{N - n},\qbs{n}) + \eb{N}{k} + \eb{N}{m}$, so with one class-$k$ and one class-$m$ customer. Note that we act as if the class-$k$ customer enters service when the high priority busy period generated by the class-$m$ customer finishes. This is feasible, since the order in which the customers are served does not alter the duration of a high priority busy period, cf. \eqref{eqn:three-class_f_3;i_2,i_1}. The remaining first passage probabilities for the case $n = 1$ are computed as, for $k \ge 2$,
\begin{equation}%
\f{k}{i_1} = 1 - \sum_{j_1 = 0}^{i_1 - 1} \g{k}{j_1}, \quad i_1 > 0. \label{eqn:N-class_f_k;i_1}
\end{equation}%

The equilibrium probabilities of the $N$-class system follow again by counting excursions as done for the two and three-class systems. The number of excursions per time unit that end in state $(\zerob{N - n},\qbs{n})$ is equal to $p(\zerob{N - n},q_n + 1,\qbs{n - 1})\mu_n$. This number is also equal to the excursions starting from class-$n$ level $q_n$ or lower per time unit that end in state $(\zerob{N - n},\qbs{n})$. A fraction $\g{n}{\ib{n - 1}}$ of the excursions starting in $(\zerob{N - n},q_n,\qbs{n - 1} - \ib{n - 1})$ by a class-$n$ arrival end in $(\zerob{N - n},\qbs{n})$. Excursions to class-$n$ levels higher than $q_n$ starting in state $(\zerob{N - n},\qbs{n} - \ib{n})$ by a class-$m, ~ m = n + 1,\ldots,N$ arrival reach, with probability $\f{m}{\ib{n}} - \g{m}{\ib{n}}$, level $q_n$ in state $(\zerob{N - n},\qbs{n})$ at first return to class-$n$ level $q_n$. We have for $n = 1,2,\ldots,N$,
\begin{align}%
&p(\zerob{N - n},q_n + 1,\qbs{n - 1})\mu_n \notag \\
&= \sum_{\ib{n} = \zerob{n}}^{\qbs{n}} p(\zerob{N - n},\qbs{n} - \ib{n}) \sum_{m = n + 1}^N \la_m (\f{m}{\ib{n}} - \g{m}{\ib{n}}) \notag \\
&+ \sum_{\ib{n - 1} = \zerob{n - 1}}^{\qbs{n - 1}} p(\zerob{N - n},q_n,\qbs{n - 1} - \ib{n - 1}) \la_n \g{n}{\ib{n - 1}}, \quad \qbs{n} \ge \zerob{n}, \label{eqn:N-class_p_n}
\end{align}%
which can be solved recursively, starting from $p(\zerob{N}) = 1 - \rho$. Note that for $n = 1$ the second term on the right-hand side of \eqref{eqn:N-class_p_n} becomes $p(\zerob{N - 1},q_1)\la_1$ and for $n = N$, the first term on the right-hand side reduces to 0.

\begin{remark}%
The above algorithm to determine the equilibrium probabilities involves subtractions in some equations, see e.g.~\eqref{eqn:N-class_f_k;i_1}, which may possibly lead to loss of significant digits and instability. However, in all experiments we observed numerically stable results.
\end{remark}%


\section{Application in spare parts logistics}%
\label{sec:application}%

Our interest in the joint queue length distribution arose from a spare parts supply problem for repairable parts sharing the same repair shop. For this problem, we apply our method, based on the matrix-analytic approach, to demonstrate the influence of assigning repair priorities on the performance of the system.

There are $M$ identical machines and each machine contains three different subsystems, numbered $1,2,3$. Each subsystem $n$ consists of $Z_n$ identical parts in parallel. We refer to the parts of subsystem $n$ as parts of Stock-Keeping Unit $n$ (SKU $n$). For each subsystem, $k_n < Z_n$ parts have to function. That is, we have redundancy, and the redundant parts are in ``cold standby''. This is called a ``$k_n$-out-of-$Z_n$'' setup. We have $k_n$ functioning parts per subsystem and only these parts are subject to failure. When one part fails, another one can immediately take over the necessary functions. An example of such a subsystem is the board computer of an airplane, where this critical component is duplicated and in an idle mode to accommodate possible failures, here, $k_n = 1$ and $Z_n = 2$. Other typical systems with this structure can be found in \cite{sherbrooke2004_system_availability}.

\begin{figure}%
\centering%
\includegraphics{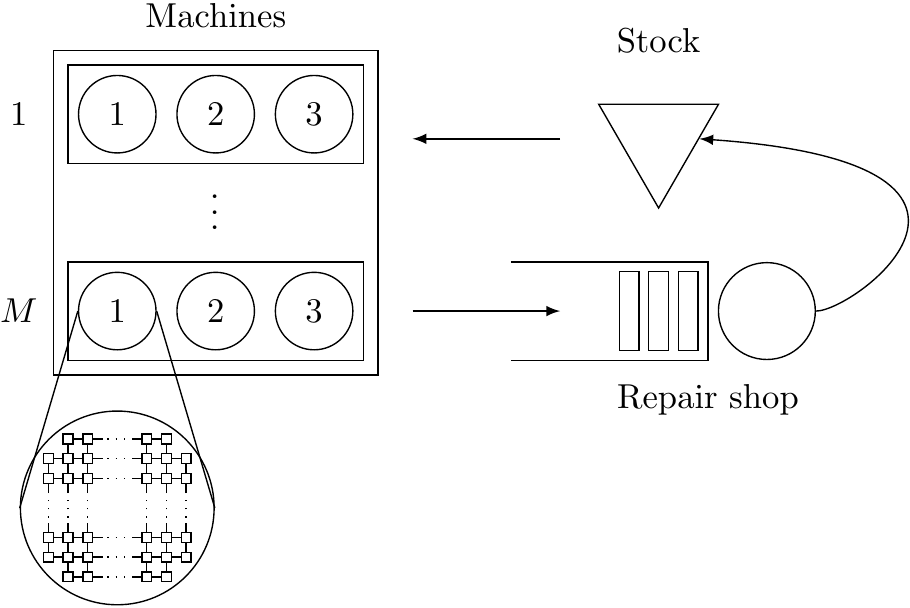}
\caption{Example of a simple spare parts supply system.}%
\label{fig:application}%
\end{figure}

A machine is only working when all three subsystems are working. When one of the functioning parts fails, a redundant part takes over its function and a service engineer takes a new part from a stock of parts and replaces the failed one. The failed part is then sent to a single server repair facility. Part and repair requests are served on a first come first serve basis. The repair time for a part of SKU $n$ is exponentially distributed with rate $\mu_n$; the delivery and replacement times are small and can be neglected. We assume that failures of parts of SKU $n$ occur according to a Poisson process with rate $\la_n$. This approximation, which is the only one needed, is valid when $M$, the total number of machines in the system, is large and when the fraction of working machines is high. After repair the broken parts are assumed to be as good as new and they are put back to stock. The stock of SKU $n$ at time instant $t = 0$ is denoted by $S_n$. We call the amount $S_n$ the basestock level for SKU $n$ parts. The system is shown in Figure~\ref{fig:application}.

Let us define the system availability as the average fraction of working machines:
\begin{equation}%
A(S_1,S_2,S_3) = \frac{1}{M} \sum_{m = 1}^M \Prob{\textup{Machine $m$ is working}}. \label{eqn:availability}
\end{equation}%
The number of backorders of SKU $n$ parts is given by $(q_n - S_n)^+$, where $(x)^+ = \max(0,x)$ and $q_n$ is the number of SKU $n$ parts in repair. Define $E_n$ as the number of `empty' spots in a given subsystem $n$ of any of the $M$ machines. Then, by conditioning on the number of parts of SKU $n$ in repair we obtain, with $s \le Z_n$,
\begin{equation}%
\Prob{E_n = s \mid q_n \text{ in repair}} = \begin{cases}%
\hfil 1, & q_n - S_n < s, ~ s = 0, \\
\hfil 0, & q_n - S_n < s, ~ s > 0, \\
\binom{Z_n}{s} \binom{Z_n(M - 1)}{q_n - S_n - s} / \binom{M Z_n}{q_n - S_n}, & q_n - S_n \ge s.
\end{cases}%
\end{equation}%
In terms of the joint queue length distribution, the system availability can be written as
\begin{align}%
A(S_1,S_2,S_3) = \sum_{q_3,q_2,q_1 \ge 0} \Bigl( \prod_{n = 1}^3  \sum_{s = 0}^{Z_n - k_n} \Prob{E_n = s \mid q_n \text{ in repair}} \Bigr) p(q_3,q_2,q_1). \label{eqn:availability_detailed}
\end{align}%
The expression \eqref{eqn:availability_detailed} determines the system availability much better than other approximations proposed in the literature; e.g., the system availability defined in \cite{sherbrooke2004_system_availability} only uses information on the mean number of backorders. The matrix-analytic method makes it possible to use the detailed distribution of the number of parts in repair. To demonstrate the approach we execute a set of experiments with the following parameters: $M = 100$, and $Z_n = 4$, $k_n = 2$, $\la_n = n/300$ and $\mu_n = (4 - n)/\be$ for $n = 1,2,3$, where $\be$ is chosen such that $\sum_{n = 1}^3 \la_n / \mu_n = \rho$.

We wish to compute the joint queue length distribution such that $\sum_{q_3,q_2,q_1} p(q_3,q_2,q_1) > 1 - \epsilon$ with $\epsilon$ a small positive number. We do this by computing the equilibrium probabilities of the states in a discrete three-dimensional cuboid $\mathcal{C}$ with states $\{0,\ldots,c_3\} \times \{0,\ldots,c_2\} \times \{0,\ldots,c_1\}$. For the sake of clarity, we briefly introduce the marginal queue length distribution of class-$n$ customers as $p_n(\cdot)$. We specify the construction of $\mathcal{C}$ in more detail. The bound $c_3$ is computed from the $M/M/1$ system with only class-$3$ customers, such that $\sum_{q_3 = 0}^{c_3} p_3(q_3) > 1 - \epsilon$, which leads to $c_3 = \lceil \frac{\log{\epsilon}}{\log{\la_3/\mu_3}} - 1 \rceil$. The bound $c_2$ is obtained through a priority system with class-$3,2$ customers such that the sum of the marginal probabilities for class-$2$ customers is very close to 1. That is, $\sum_{q_{2} = 0}^{c_{2}} p_2(q_2) > 1 - \epsilon$. Conveniently, the marginal queue length distribution $p_2(\cdot)$ can be derived directly from the joint equilibrium probabilities $p(0,\cdot)$ of the priority queueing system with class-$3,2$ customers via the relation $\la_2 p_2(q_2 - 1) = \mu_2 p(0,q_2)$. Thus, this allows us to estimate the bound $c_2$ without having to compute all joint equilibrium probabilities $p(q_3,q_2)$. The final bound $c_1$ can be found iteratively until $\sum_{q_3,q_2,q_1} p(q_3,q_2,q_1) > 1 - \epsilon$ or using the same method as for the bound $c_2$. Naturally, this method of constructing $\mathcal{C}$ extends to an arbitrary number of classes.

In Table~\ref{tbl:application_results} we list the system availability according to \eqref{eqn:availability_detailed} for different utilization rates of the repair shop and different priority assignments. The basestock levels $S_n$ depend on the mean queue lengths, i.e., we set $S_n = \lfloor \E{Q_n} \rfloor, ~ n = 1,2,3$, where $Q_n$ is the queue length of SKU $n$ parts. The algorithm for the 3-class system was executed using Java 8.0 on a PC with an Intel Core i7-3770 CPU and 16 GB RAM. The computation times mentioned in Table~\ref{tbl:application_results} depend on the number of states with significant probability mass, i.e., on the load of the system, the priority assignment and, naturally, the parameter value of $\epsilon$. For these experiments, we have selected $\epsilon = 10^{-6}$.

\begin{table}%
\centering%
\begin{tabular}{|r|ccc|rrr|c|r|}%
\hline%
Util. & \multicolumn{3}{c}{Priorities} & \multicolumn{3}{|c|}{Mean queue length} & Avail. & Comp. \\%
$\rho$ & $r_1$ & $r_2$ & $r_3$ & SKU 1 & SKU 2 & SKU 3 & $A$ & time (s)\\%
\hline \hline%
0.90 & H & M & L &  0.0744 &  0.3015 &  7.3244 & 0.9999 &  0.07 \\ 
     & H & L & M &  0.0744 & 10.7998 &  2.0752 & 0.9996 &  0.17 \\ 
     & M & H & L &  0.1333 &  0.2621 &  7.3244 & 0.9999 &  0.02 \\ 
     & M & L & H &  1.3408 & 10.7998 &  1.6531 & 0.9996 &  0.78 \\ 
     & L & H & M &  9.6132 &  0.2621 &  4.1643 & 0.9995 &  0.60 \\ 
     & L & M & H &  9.6132 &  5.2846 &  1.6531 & 0.9994 &  4.28 \\ 
\hline%
0.95 & H & M & L &  0.0788 &  0.3261 & 15.6437 & 0.9995 &  0.17 \\ 
     & H & L & M &  0.0788 & 26.5995 &  2.5070 & 0.9965 &  2.03 \\ 
     & M & H & L &  0.1467 &  0.2808 & 15.6437 & 0.9995 &  0.04 \\ 
     & M & L & H &  1.8359 & 26.5995 &  1.9213 & 0.9965 &  7.50 \\ 
     & L & H & M & 28.7923 &  0.2808 &  6.0938 & 0.9930 & 11.02 \\ 
     & L & M & H & 28.7923 &  8.6257 &  1.9213 & 0.9928 & 34.35 \\ 
\hline%
\end{tabular}%
\caption{System availability for different combinations of repair shop utilizations and priority assignments. We use $\epsilon = 10^{-6}$. The variable $r_n$ indicates the priority of SKU $n$ parts, either high (H), medium (M) or low (L).}%
\label{tbl:application_results}%
\end{table}%

Table~\ref{tbl:application_results} shows that we have a fast numerical method to compute the availability for different priority assignments and particular choices of the basestock levels. This method can easily be exploited in a procedure to optimize the priority assignment and basestock level; e.g.~in order to maximize system availability under a given budget for spare parts (cf.~\cite{adan2009_spare_parts_problem} which considers a slightly different setting with equal repair rates for all SKU's).


\section{Conclusion and extensions}%
\label{sec:conclusion}%

We have developed for the $M/M/1$ preemptive priority system with $N$ customer classes and class-dependent service rates a method for the exact determination of the joint equilibrium queue length distribution. This method is based on the matrix-analytic method as the embedded Markov processes are of the $M/G/1$ type. Key to this approach are first passage probabilities, computed by one-step analysis.

We applied the exact solution method to a spare parts logistics problem where repairable parts share the same repair shop, and showed that this method produces accurate results in the order of seconds.

\begin{figure}%
\centering%
\includegraphics{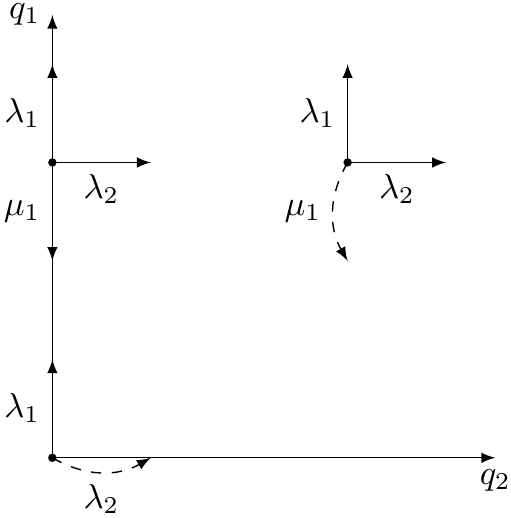}
\caption{Transition rate diagram of the non-preemptive $M/M/1$ priority system of the states $(0,q_1,1)$ with $q_1 > 0$ and including state $(0,0,0)$. The dashed arrows indicate a transition to a state with $s = 2$, i.e.~a state with a class-2 customer in service.}%
\label{fig:non-preemptive_trd_s_1}%
\end{figure}

We next sketch how the method can be extended to an $M/M/1$ non-preemptive priority system. In the non-preemptive case one identifies the customer currently in service by adding another variable to the state description. For the two-class system, the state description becomes $(q_2,q_1,s)$ where $s \in \{1,2\}$ indicates the class of the customer in service and $s = 0$ indicates no customer in service. By defining class-2 level $q_2$ as the set of states with $q_2$ class-2 customers, one can again count the number of excursions per time unit that start from class-2 level $q_2$ and reach levels higher than $q_2$ to finally end at state $(q_2,q_1,2)$. The states with a class-1 customer in service can only be reached from the states $(q_2,q_1,1)$ or $(0,0,0)$ and thus the equilibrium probabilities of these states can be recursively determined for $q_2 > 0$ immediately from the boundary probabilities of class-2 level 0, see Figure~\ref{fig:non-preemptive_trd_s_1}. One finds the equilibrium probabilities of class-2 level 0, starting from $p(0,0,0) = 1 - \rho$, by embedding the Markov process on class-2 level 0 and again counting excursions. Notice that the approach is very similar to the one for the preemptive case and only requires the computation of equilibrium probabilities of the states $(q_2,q_1,1)$ as an additional step.

\bibliographystyle{plain}%
\bibliography{bib-priorities}%

\end{document}